\documentclass[12pt,draft]{extartic}

\usepackage{amsmath,amsthm,amssymb,calc}
\usepackage[dvips]{graphics, graphicx}

\usepackage[english]{babel}                         %eng

\newcounter{num}[section]

\newcommand{\Th}{\refstepcounter{num}
{\bf Theorem \arabic{section}.\arabic{num} }}

\newcommand{\Lemma}{\refstepcounter{num}
{\bf Lemma \arabic{section}.\arabic{num} }}

\newcommand{\Pred}{\refstepcounter{num}
{\bf Proposition \arabic{section}.\arabic{num} }}

\newcommand{\Cor}{\refstepcounter{num}
{\bf Corollary \arabic{section}.\arabic{num} }}

\newcommand{\Note}{\refstepcounter{num}
{\it Note \arabic{section}.\arabic{num} }}

\newcommand{\Exm}{\refstepcounter{num}
{\bf Example \arabic{section}.\arabic{num} }}

\newcommand{\St}{\refstepcounter{num}
{\bf Statement \arabic{section}.\arabic{num} }}

\newcommand{\Def}{\refstepcounter{num}
{\it Definition \arabic{section}.\arabic{num} }}

\newcommand{\Proof}{{\bf Proof. }}

\def\eps{\varepsilon}
\def\_phi{\varphi}

\def\a{\alpha}
\def\d{\delta}
\def\la{\lambda}

\def\F{\widehat}
\def\L{\Lambda}
\def\m{\times}

\def\ov{\overline}

\def\C{{\mathbb C}}

\def\Z_N{{\mathbb Z}_N}
\def\Z{{\mathbb Z}}

\def\Span{{\rm Span\,}}

\def\Gr{{\mathbf G}}

\def\dim{{\rm dim }}

\def\l{\left}
\def\r{\right}

\def\Spec{{\rm Spec\,}}
\def\oC{{\rm C}}
\def\oI{{\rm I}}
\def\oM{{\rm M}}
\def\oP{{\rm P}}
\def\oT{{\rm T}}
\def\oS{{\rm S}}
\def\supp{{\rm supp\,}}
\def\tr{{\rm tr\,}}

\author{Shkredov I.D.}
\title{Some applications of W. Rudin's inequality to problems of combinatorial number theory
\footnote{ This work was supported by RFFI grant no. 06-01-00383,
P. Deligne's grant (Balsan's fond 2004),
President's of Russian Federation grant N MK--1959.2009.1
and grant of Leading Scientific Schools 691.2008.1}}
\date{}
\begin{document}
\maketitle

\begin{center}
 Annotation.
\end{center}

{\it \small
    In the paper
    we obtain
    some new applications of well--known W. Rudin's theorem concerning lacunary
    series to problems of combinatorial number theory.
    %were obtained.
    We generalize a result of M.--C. Chang on $L_2 (\L)$--norm
    of
    Fourier coefficients
    %of the characteristic function
    of a set
%    with
%    and
    (here $\L$ is a dissociated set),
    and prove a dual version of the theorem.
    Our
    %basic method of investigation
    main instrument
    is computing of eigenvalues of some operators.
}
\\
%\\
%\\

\refstepcounter{section}
\label{introduction}

\begin{center}
{\bf \arabic{section}. Introduction.}
\end{center}

A well--known theorem of W. Rudin concerning lacunary series (see \cite{Rudin_book,Rudin}) states
%that
for any complex function with condition that support of its Fourier transform
belongs to a set "without arithmetical structure"\,
all $L_p$---norms, $p\ge 2$ are equivalent to $L_2$---norm.
Let us formulate the last result more precisely.
Let $\Gr=(\Gr,+)$ be an
%finite
 Abelian group with additive group operation $+$.
By $\F{\Gr}$ denote the Pontryagin dual of $\Gr$.
In other words $\F{\Gr}$ is the group of homomorphisms $\xi$ from $\Gr$ to $\mathbf{R} / \mathbf{Z}$,
$\xi : x \to \xi \cdot x$.
A set $\L \subseteq \F{\Gr}$, $\L = \{ \la_1,\dots, \la_{|\L|} \}$
is called {\it dissociated} if any identity of the form
$$
    \sum_{i=1}^{|\L|} \eps_i \la_i = 0 \,,
$$
where $\eps_i \in \{ 0, \pm 1 \}$ implies that all $\eps_i$ are equal zero.

\Th
   \label{t:Rudin}
   {\it
        Let $\Gr$ be a finite Abelian group.
        There exists an absolute constant $C>0$
        such that for any dissociated set $\L \subseteq \Gr$,
        any complex numbers $a_n \in \C$, and all positive integers $p \ge 2$
        the following inequality holds
        \begin{equation}\label{f:Rudin}
            \frac{1}{|\Gr|} \sum_{x \in \Gr} \left| \sum_{\xi\in \L} a_\xi e^{2\pi i \xi \cdot x} \right|^p \le
                (C \sqrt{p})^p \left( \sum_{\xi\in \L} |a_\xi|^2 \right)^{p/2} \,.
        \end{equation}
   }

Theorem \ref{t:Rudin} is widely used in combinatorial number theory,
mainly in solution of so--called "inverse"\, problems (see books \cite{Nathanson,Tao_Vu_book}).
Thus, using Rudin's theorem M.--C. Chang proved the following result (see \cite{Ch_Fr,Green_Chang2,Ruzsa_Green_Fr}).

\Th
{\it
    Let $\Gr$ be a finite Abelian group, and $\d \in (0,1]$ be a real number.
    Let also
    $\L \subseteq \Gr$ be a dissociated set, $S\subseteq \Gr$ be an arbitrary set, $|S| = \d |\Gr|$.
    Then
    \begin{equation}\label{F:Chang_l_2}
        \sum_{\xi \in \L} \l| \sum_{x\in S} e^{2\pi i \xi \cdot x} \r|^2 \le C_1 |S|^2 \log (1/\d) \,,
    \end{equation}
    where $C_1 > 0$ is an absolute constant.
}
\label{t:Chang}

Theorem \ref{t:Chang} was repeatedly used in problems concerning
large exponential sums
(see papers \cite{Ch_Fr,GreenA+A,Ruzsa_Green_Fr,Sanders_L1_Chang,Sanders_Fr_new} and others).
On sets of large exponential sums see \cite{Gow_surv}, \cite{Green_Chang_exact}--\cite{Green_Chang2}, \cite{Tao_Vu_book},
\cite{Freiman_Yudin}--\cite{Ball}, \cite{Sh_dokl_exp1}--\cite{Shkr_Bourgain_group}.
In paper  \cite{BourgainA+A} Rudin's theorem
was used for finding arithmetic progressions in sumsets,
and for studying properties of sets with small doubling
(see \cite{Sanders_m}, \cite{Sh_doubling}).
In the paper we obtain new applications of the beautiful theorem
to problems of combinatorial number theory.
Firstly, we prove a sharp generalization of Chang's theorem.

\Th
{\it
    Let $\Gr$ be a finite Abelian group, $\d \in (0,1]$ be a real number,
    and $l$ be a positive integer, $l\ge 2$.
    Suppose that $\L \subseteq \Gr$ is a dissociated set, and $S\subseteq \Gr$ is an arbitrary set, $|S| = \d |\Gr|$.
    Then
\begin{equation}\label{f:my_Chang_gen}
    \sum_{\xi\in \L} \l| \sum_{x\in S} e^{2\pi i \xi \cdot x} \r|^{l+1}
        \le
            C_2 |S| \l( \sum_{\xi \neq 0} \l| \sum_{x\in S} e^{2\pi i \xi \cdot x} \r|^{2l} \r)^{1/2} \cdot \log^{1/2} (1/\d) \,,
\end{equation}
    where $C_2 > 0$ is an absolute constant.
}
\label{t:my_Chang_gen}

Secondly, we obtain a "dual"\, version of Theorem \ref{t:my_Chang_gen}.
%(дуальность операции замены преобразования Фурье на свертку обсуждалась in книге \cite{Tao_Vu_book}).

\Th
{\it
    Let $\Gr$ be a finite Abelian group, and $l$ be  a positive integer, $l\ge 2$.
    Suppose that $\L \subseteq \Gr$ is a dissociated set, and $S_1,\dots, S_l \subseteq \Gr$ are arbitrary sets.
    Then
\begin{equation}\label{f:my_Chang_dual}
    \sum_{x\in \L} (S_1 * S_2 * \dots * S_l)^2 (x)
        \le
                C_3
            \frac{|S_l|}{N} \l( \sum_{\xi} \prod_{j=1}^{l-1} \l| \sum_{x\in S_j} e^{2\pi i \xi \cdot x} \r|^2 \r) \cdot \log N \,,
\end{equation}
    where $C_3 > 0$ is an absolute constant.
}
\label{t:my_Chang_dual}

Structure of the paper is the following.
In section \ref{operators_sec} we introduce operators (matrices) $\oT^{\_phi}_\psi$ and $\oS^{\_phi}_\psi$.
This family of matrices is our main instrument of investigation.
We study spectrums of matrices $\oT^{\_phi}_\psi$ and $\oS^{\_phi}_\psi$
and its eigenvectors.
In section  \ref{Chang_gen} we reformulated Rudin's theorem on the language of eigenvalues
of operators $\oT^{\_phi}_\L$, where $\L$ is a dissociated set and $\_phi$ is a function
(see Proposition \ref{st:mu_1_L_S})
and derive Theorems \ref{t:my_Chang_gen}, \ref{t:my_Chang_dual}.
Also we
%demonstrate
show
that Theorem \ref{t:my_Chang_gen} is sharp and
for some choice of parameters
Theorem \ref{t:my_Chang_dual} is also sharp.

Let us say a few words about the notation.
If $S\subseteq \Gr$ is a set then we will write $S(x)$ for the characteristic function.
In other words $S(x) = 1$ if $x\in S$ and zero otherwise.
By $\log$ denote logarithm base two.
Sings $\ll$ and $\gg$ are usual Vinogradov's symbols.
If $n$  is a positive integer then we will write $[n]$ for the segment $\{ 1,2,\dots,n \}$.

%Автор выражает благодарность Математическому Исследовательскому Институту (MSRI)
%за его гостеприимство and создание прекрасных условий для работы.
%Также хотелось бы поблагодарить доктора физико--математических наук,
%профессора Н.\,Г.\,Мощевитина and доктора физико--математических наук, профессора С.\,В.\,Конягина
%за внимание к работе and полезные обсуждения.
The author is grateful to N.\,G.\,Moshchevitin and S.\,V.\,Konyagin
for their attention to the work.
Also he would like to thanks S. Yekhanin and Kunal Talwar for useful discussions.

\refstepcounter{section}
\label{operators_sec}

\begin{center}
{\bf \arabic{section}. Operators $\oT^{\_phi}_\psi$ and $\oS^{\_phi}_\psi$.}
\end{center}

%%В дальнейшем изложении
%%%Нам понадобится анализ Фурье.
We need in Fourier analysis in our proof.
%%%%Recall some notation
%%%%Напомним некоторые сведения из анализа Фурье.
%%Для конечной абелевой группы $\Gr=(\Gr,+)$ с аддитивной групповой операцией $+$
%%обозначим через $\F{\Gr}$ двойственную группу для $\Gr$.
%% In other words let $\F{\Gr}$ --- group гомоморфизмов $\xi$ из $\Gr$ in $\mathbf{R} / \mathbf{Z}$, $\xi : x \to \xi \cdot x$.
Let $\Gr$ be a finite Abelian group, $N:=|\Gr|$.
It is well--known that
in the case the dual group $\F{\Gr}$ is isomarphic to  $\Gr$.
Let also $f$ be a function from $\Gr$ to $\C$.
% с конечным носителем.
By $(\Phi f) (\xi) = \F{f}(\xi)$ denote the Fourier transformation of $f$
\begin{equation}\label{F:Fourier}
  (\Phi f) (\xi) = \F{f}(\xi) =  \sum_{x \in \Gr} f(x) e( -\xi \cdot x) \,,
\end{equation}
where $e(x) = e^{2\pi i x}$.
We will use the following basic facts
\begin{equation}\label{F_Par}
    \| f \|^2_2 := \sum_{x\in \Gr} |f(x)|^2
        =
            \frac{1}{N} \sum_{\xi \in \F{\Gr}} |\widehat{f} (\xi)|^2
                =
                    \frac{1}{N} \| \F{f} \|^2_2 \,.
\end{equation}
\begin{equation}\label{F_Par_sc}
    \langle f,g \rangle
        :=
            \sum_{x\in \Gr} f(x) \overline{g(x)}
                = \frac{1}{N} \sum_{\xi \in \F{\Gr}} \F{f}(\xi) \overline{\F{g}(\xi)}
                    = \frac{1}{N} \langle \F{f}, \F{g} \rangle \,.
\end{equation}
\begin{equation}\label{svertka}
    \sum_{y\in \Gr} |\sum_{x\in \Gr} f(x) g(y-x) |^2
        = \frac{1}{N} \sum_{\xi \in \F{\Gr}} |\widehat{f} (\xi)|^2 |\widehat{g} (\xi)|^2 \,.
\end{equation}
\begin{equation}\label{f:inverse}
    f(x) = \frac{1}{N} \sum_{\xi \in \F{\Gr}} \F{f}(\xi) e(\xi \cdot x) \,.
\end{equation}
%\begin{equation}\label{f:U_2}
%    \| f \|_{U_2} = \l( \frac{1}{N^3} \sum_{x,y,z} f(x) \ov{f(x+y)}\, \ov{f(x+z)} f(x+y+z)\r)^{1/4}
%                        = \left( \frac{1}{N^4} \sum_{r\in \Gr} |\F{f} (r)|^4 \right)^{1/4} \,.
%\end{equation}
If
$$
    (f*g) (x) := \sum_{y\in \Gr} f(y) g(x-y)
$$
 then
\begin{equation}\label{f:F_svertka}
    \F{f*g} = \F{f} \F{g} \quad \mbox{ and } \quad (\F{fg}) (x) = \frac{1}{N} (\F{f} * \F{g}) (x) \,.
\end{equation}
By $*_{l-1}$ denote the result of using convolution $l$ times.
Let also $(f *_{0} f) (x) := f(x)$.
% ЭТА свертка переводится преобр. Фурье просто in произведение f на g.

We will write
$\sum_s$ instead of $\sum_{s\in \Gr}$
 and
$\sum_{\xi}$ instead of $\sum_{\xi\in \F{\Gr}}$
for brevity.

\Def Let $\_phi,\psi : \Gr \to \C$ be two functions.
By $\oT^{\_phi}_\psi$ denote the following operator on the space of functions $\Gr^{\C}$
\begin{equation}\label{F:T}
    (\oT^{\_phi}_\psi f ) (x) = \psi(x) (\F{\_phi^c} * f) (x) = \psi(x) \F{ \_phi^c \F{f^c}} (x) \,,
\end{equation}
where $f$ is an arbitrary complex function on $\Gr$ and
$f^c$ is the function $f^c (x) = f(-x)$.
Also we need in "more symmetric"\, (see identity (\ref{F:S^*}) below) operator $\oS^{\_phi}_\psi$
\begin{equation}\label{F:S}
    (\oS^{\_phi}_\psi f ) (x) = \psi(x) (\F{\_phi^c} * \ov{\psi} f) (x) = \psi(x) (\_phi^c \F{\ov{\psi^c} f^c})\,\,\F{}\, (x) \,.
\end{equation}

In particular if $\psi \equiv 1$ then $\oT^{\_phi}_\psi$ if the convolution operator
and if $\_phi \equiv 1$ then $\oT^{\_phi}_\psi$ is the operator of multiplication by the function.

Let us express the operators $\oT^{\_phi}_\psi$, $\oS^{\_phi}_\psi$
as composition of more simple operators.
Let
$$
    (\oC f) (x) = f(-x) = f^c (x)
$$
and for any complex function $\rho : \Gr \to \C$ let
$$
    (\oP_\rho f) (x) = \rho (x) f(x) \,.
$$
Clearly, $\oC^2 = \oI$ is the identity operator and for any two functions $\rho_1$, $\rho_2$
the following holds
$\oP_{\rho_1} \oP_{\rho_2} = \oP_{\rho_1 \rho_2}$.
We have
\begin{equation}\label{F:operators_identities}
    \oC \Phi = \Phi \oC \,, \quad \oC \oP_\rho = \oP_{\rho^c} \oC \,, \quad \Phi^2 = N \cdot \oC \,.
\end{equation}
%Легко видеть, что
The last formulas imply
\begin{equation}\label{F:T_simple_op}
    \oT^{\_phi}_\psi = \oP_{\psi} \Phi \oP_{\_phi^c} \Phi \oC = \oC \oP_{\psi^c} \Phi \oP_{\_phi} \Phi
        =
            \oP_{\psi} \Phi \oC \oP_{\_phi} \Phi
\end{equation}
and
\begin{equation}\label{F:S_simple_op}
    \oS^{\_phi}_\psi
        = \oP_{\psi} \Phi \oP_{\_phi^c} \Phi \oP_{\ov{\psi}^c} \oC
            = \oC \oP_{\psi^c} \Phi \oP_{\_phi} \Phi \oP_{\ov{\psi}}
                =
                    \oP_{\psi} \Phi \oC \oP_{\_phi} \Phi \oP_{\ov{\psi}} \,.
\end{equation}
We will need in more formulas.
%Отметим еще несколько формул, которые пригодятся нам in дальнейшем.
We have
\begin{equation}\label{F:clever_convolution}
    \oC (a*b) = \oC a * \oC b \,, \quad \langle a, \ov{b}*c \rangle = \langle b, \ov{a}*\oC c \rangle
\end{equation}
for any functions $a,b$ and $c$.
Let $\mathcal{C}$ be the conjugation operator.
Clearly, $\mathcal{C}^2 = \oI$.
% для любой function $f : \Gr \to \C$.
Besides
\begin{equation}\label{F:complex_op}
    \mathcal{C} \oC = \oC \mathcal{C} \,, \quad \mathcal{C} \Phi = \Phi \mathcal{C} \oC = \Phi \oC \mathcal{C} \,, \quad
    \Phi \mathcal{C} = \mathcal{C} \Phi \oC = \mathcal{C} \oC \Phi \,.
\end{equation}

Now let us find operators $(\oT^{\_phi}_\psi)^*$ and $(\oS^{\_phi}_\psi)^*$.
First of all, note that formula (\ref{F_Par_sc}) is equivalent to identity
$\Phi^* = \Phi \oC = \oC \Phi$.
Secondly, we have $\oP^*_{\rho} = \oP_{\ov{\rho}}$ and $\oC^* = \oC$.
From this and identities (\ref{F:operators_identities}), we get
\begin{equation}\label{F:T^*}
    (\oT^{\_phi}_\psi)^* = \oC \Phi \oP_{\ov{\_phi}} \Phi \oP_{\ov{\psi}} = \Phi \oP_{\ov{\_phi^c}} \Phi \oP_{\ov{\psi^c}} \oC \,,
\end{equation}
and
\begin{equation}\label{F:S^*}
    (\oS^{\_phi}_\psi)^*
        = \oC \oP_{\psi^c} \Phi \oP_{\ov{\_phi}} \Phi \oP_{\ov{\psi}}
            = \oP_{\psi} \Phi \oP_{\ov{\_phi^c}} \Phi \oP_{\ov{\psi^c}} \oC
                = \oS^{\ov{\_phi}}_\psi \,.
\end{equation}
In particular an operator $\oS^{\_phi}_\psi$ is hermitian, provided by $\_phi$ is a real function.
Finally, note that
\begin{equation}\label{F:T_action}
    \langle \oT^{\_phi}_\psi u, v \rangle = \sum_x \psi(x) \ov{v} (x) (\F{\_phi^c} * u) (x)
        =
            \sum_x \_phi(x) \F{u}(x) \ov{ (\ov{\psi} v)\,\,\F{}\, (x) }
\end{equation}
and
\begin{equation}\label{F:S_action}
    \langle \oS^{\_phi}_\psi u, v \rangle = \sum_x \psi(x) \ov{v} (x) (\F{\_phi^c} * \ov{\psi} u) (x)
        =
            \sum_x \_phi(x) (\ov{\psi} u)\,\,\F{}\, (x) \ov{ (\ov{\psi} v)\,\,\F{}\, (x) } \,.
\end{equation}
To obtain the last formulas we have used identity (\ref{F_Par_sc}).

We need in some statements from linear algebra.
Let $M$ be an arbitrary matrix $n\m n$, $M=(m_{ij})_{i,j=1}^n$.
By $\Spec (M)$ (or $\Spec M$) denote the multiset of $n$ eigenvalues of $M$.
For example $\Spec \oP_\rho = \{ \rho (x_1), \dots, \rho (x_N) \}$, where $\{ x_j \}_{j=1}^N = \Gr$
and
\begin{displaymath}
\label{F:d_c}
  \d_c (x) =
  \left\{ \begin{array}{ll}
        1, &  \mbox{ if } x=c \\
        0, &  \mbox{ otherwise. }
  \end{array} \right.
\end{displaymath}
are orthonormal eigenfunctions of operator $\oP_\rho$.
Lemma below can be find in \cite{La} p. 104.

\Lemma
{\it
    For any matrices $M_1$ and $M_2$ (possibly singular), we have
\begin{equation}\label{f:comm_spec}
    \Spec (M_1 M_2) = \Spec (M_2 M_1) \,.
\end{equation}
}
\label{l:comm_matrix}
%
%
%
%%%The equality in (\ref{f:comm_spec}) in the sense of multisets.

Denote by $\mu_j (M)$, $j=1,\dots,n$ the eigenvalues of $M$.
Let us remind that a matrix $M$ is called {\it symmetric} (or hermitian) if $M^* = M$,
where $M^* = (m^*_{ij})_{i,j=1}^n$, $m^*_{ij} = \ov{m}_{ji}$.
It is well--known that all eigenvalues of such matrices are real.
Suppose that the eigenvalues are arranged in order of magnitude
\begin{equation}\label{F:mu_ordering}
    \mu_1 (M) \ge \mu_2 (M) \ge \dots \ge \mu_n (M) \,.
\end{equation}
Courant--Fischer's theorem   (see e.g. \cite{La}, p. 115) gives us
the characterization of eigenvalues of symmetric matrices.

\Th
{\it
    Let $M$ be a symmetric matrix $n\m n$.
    Suppose that its eigenvalues $\mu_j$ are arranged like in (\ref{F:mu_ordering}).
    Then for any $r=0,1,\dots,n-1$, we have
    $$
        \mu_{r+1} (M) = \min_{v_1,\dots,v_r}\, \max_{f \,:\, \langle f,v_j \rangle = 0,\, \| f \|_2 = 1} \, \langle M f,f \rangle \,,
    $$
    where the minimum is taken over all systems of orthogonal vectors $v_1,\dots,v_r$.
    Besides, for any  vectors $u$ and $v$ the following holds
    $$
        | \langle M u, v \rangle |
            \le
                |\mu_1 (M)| \cdot \| u \|_2 \| v\|_2 \,.
    $$
}
\label{t:KF_matrix}

The
%situation
case
$\Gr = \Z / p\Z := \Z_p$, where $p$ is a prime number
is important for applications.
Below, we will obtain some statements for the situation.
In the proof we will need in Chebotar\"{e}v's theorem  (see e.g. \cite{Chebotarev}).
%Theorem Чеботарева ниже --- общеизвестна (see e.g. \cite{Chebotarev}).

\Th
{\it
    Let $\Gr = \Z_p$, and $p$ be a prime number.
    Then any submatrix  of the operator $\Phi$ is nonsingular.
}
\label{t:Chebotarev}

Theorem \ref{t:Chebotarev} implies the following result (see e.g. \cite{Tao_Vu_book} or \cite{Tao_uncertainty}).

\Th
{\it
    Let $\Gr = \Z_p$, and $p$ be a prime number.
    Let also $f:\Gr \to \C$  be an arbitrary non--zero complex function.
    Then
    \begin{equation}\label{}
        |\supp f| + |\supp \F{f}| \ge p+1 \,.
    \end{equation}
}
\label{t:Tao}

By $\Span  \{ f_j \} $ denote the linear hull of a system of vectors  $\{ f_j \}$.

\Cor
{\it
    Let $\Gr = \Z_p$, and  $p$ be a prime number.
    Suppose that  $S = \{ s_1, \dots, s_{|S|} \} \subseteq \Z_p$ is an arbitrary nonempty set.
    Let also $l$ be a positive integer, $\{ f_j \}$ be a family of functions,
    $\supp f_j \subseteq S$, and $\dim (\Span \{ f_j \}) = l$.
    Then there exist $l$ linearly independent functions $g_i \in \Span \{ f_j \}$, $i\in [l]$ such that
    $|\supp g_i| \ge p-|S|+l$.
    If $l=|S|$ then $g_i = \d_{s_i}$, $i\in [|S|]$.
}
\label{cor:Tao_th}
\\
\Proof
Choosing an appropriate basis, we find  linearly independent functions $g_i \in \Span \{ f_j \}$, $i\in [l]$ such that
$|\supp g_i| \le |S| - l + 1$,
and if $l=|S|$ then $g_i = \d_{s_i}$, $i\in [|S|]$.
Using theorem \ref{t:Tao}, we get  $|\supp g_i| \ge p-|S|+l$, $i\in [l]$.
This completes the proof.

\Pred
{\it
    We have \\
$1)~$ $\Spec (\oT^{\psi}_{\_phi}) = \Spec (\oT^{\_phi}_{\psi^c}) = \Spec (\oT^{\_phi^c}_{\psi})$. \\
$2)~$ $\Spec (\oT^{\_phi}_\psi (\oT^{\_phi}_\psi)^*) = N \cdot \Spec( \oT^{|\_phi|^2}_{|\psi|^2} )$
 and $\oT^{\_phi}_\psi (\oT^{\_phi}_\psi)^* = N \cdot \oS^{|\_phi|^2}_\psi$. \\
$3)~$ $\Spec (\oS^{|\_phi|^2}_\psi) = \Spec (\oS^{|\psi^c|^2}_\_phi) = \Spec (\oS^{|\psi|^2}_{\_phi^c})$. \\
$4)~$ $\Spec (\oS^{\_phi}_\psi) = \Spec (\oT^{\_phi}_{|\psi|^2})$.
}
\label{pred:Spec_T}
\\
\Proof
%We have
By (\ref{F:T_simple_op}), we have
\begin{equation}\label{tmp:20.06.09_1}
    \oT^{\psi}_{\_phi} = (\oP_{\_phi} \Phi) (\oP_{\psi^c} \Phi \oC)
\end{equation}
 and
$$
    (\oP_{\psi^c} \Phi \oC) (\oP_{\_phi} \Phi) = \oP_{\psi^c} \Phi \oP_{\_phi^c} \Phi \oC = \oT^{\_phi}_{\psi^c} \,.
$$
By Lemma \ref{l:comm_matrix}, we get $\Spec (\oT^{\psi}_{\_phi}) = \Spec (\oT^{\_phi}_{\psi^c})$.
Similarly, using the second identity from (\ref{F:T_simple_op}) instead of (\ref{tmp:20.06.09_1}),
we obtain $\Spec (\oT^{\psi}_{\_phi}) = \Spec (\oT^{\_phi^c}_{\psi})$.

Let us prove the second part of our proposition.
From (\ref{F:operators_identities}), (\ref{F:T_simple_op}) and (\ref{F:T^*}), we have
$$
    \oT^{\_phi}_\psi (\oT^{\_phi}_\psi)^*
        =
            \oP_{\psi} \Phi \oP_{\_phi^c} \Phi \oC
                \Phi \oP_{\ov{\_phi^c}} \Phi \oP_{\ov{\psi^c}} \oC
                    =
                        N \cdot ( \oP_{\psi} \Phi ) ( \oP_{|\_phi^c|^2} \Phi \oP_{\ov{\psi^c}} \oC )
                            = N \cdot \oS^{|\_phi|^2}_\psi \,.
$$
Further
$$
    ( \oP_{|\_phi^c|^2} \Phi \oP_{\ov{\psi^c}} \oC ) ( \oP_{\psi} \Phi )
        =
            \oP_{|\_phi^c|^2} \Phi \oP_{|\psi^c|^2 } \Phi \oC
$$
and by Lemma \ref{l:comm_matrix} again and also by the first part, we get
$$
    \Spec(\oT^{\_phi}_\psi (\oT^{\_phi}_\psi)^*)
        =
            N \cdot \Spec (\oT^{|\psi|^2}_{|\_phi^c|^2})
                =
                    N \cdot \Spec (\oT^{|\_phi|^2}_{|\psi|^2}) \,.
$$
The third part is a corollary of $1)$ and $2)$.

We need to check the last part of proposition \ref{pred:Spec_T}.
We have
$$
    \oS^{\_phi}_\psi
        = (\oP_{\psi} \Phi) (\oP_{\_phi^c} \Phi \oP_{\ov{\psi}^c} \oC) \,.
$$
Hence the multiset  $\Spec(\oS^{\_phi}_\psi)$ is coincident  with the spectrum of operator
$$
    \oP_{\_phi^c} \Phi \oP_{\ov{\psi}^c} \oC \oP_{\psi} \Phi = \oP_{\_phi^c} \Phi \oP_{|\psi^c|^2} \Phi \oC = \oT^{|\psi|^2}_{\_phi^c} \,.
$$
Using $1)$, we get  $\Spec (\oS^{\_phi}_\psi) = \Spec (\oT^{\_phi}_{|\psi|^2})$.
This completes the proof.

Note that
$\Spec (\oS^{\_phi}_\psi) \neq \Spec (\oS^{\psi}_\_phi)$
in general.

Now let us consider a class of operators $\oT^{\_phi}_\psi$ of special type.
Let $S=\{ s_1,\dots, s_{|S|} \} \subseteq \Gr$ be an arbitrary set.
We are interested in the class of operators $\oT^{\_phi}_\psi$ of the form $\oT^{\_phi}_S$,
where $\psi(x) = S(x)$ is the characteristic function of our set $S$.

Let $L(S) = \{ f ~:~ \supp f \subseteq S \}$ and $L(\ov{S}) = \{ f ~:~ \supp f \subseteq \Gr\setminus S \}$.
Clearly, the linear space
%functions
$f: \Gr \to \C$ is a direct sum of subspaces $L(S)$ and $L(\ov{S})$.
Obviously, $\dim L(S) = |S|$, $\dim L(\ov{S}) = N - |S|$ and $\oT^{\_phi}_S ( L(S) ) \subseteq L(S)$.

\Def By  $\ov{\oT}^{\_phi}_S$ denote the restriction of operator $\oT^{\_phi}_S$ onto the space $L(S)$.
Let also $\ov{\oS}^{\_phi}_S$ be the restriction of $\oS^{\_phi}_S$ onto  $L(S)$.

Clearly,$\ov{\oT}^{\_phi}_S = \ov{\oS}^{\_phi}_S$ but
$\oT^{\_phi}_S \neq \oS^{\_phi}_S$ in general.
Nevertheless, we will prove (see proposition below) that
$\Spec \oT^{\_phi}_S = \Spec \oS^{\_phi}_S$.

The matrix of operator $\ov{\oT}^{\_phi}_S$ is $( \ov{\oT}^{\_phi}_S )_{ij} = \F{\_phi} (s_i-s_j)$.
It is easy to see then $\ov{\oT}^{\_phi}_S$ is a symmetric operator, provided by $\_phi$ is a real function.

\Pred
{\it
    We have \\
$1)~$ $\Spec (\oT^{\_phi}_{S}) =
%        \Spec (\oT^{\_phi}_{S} \oP_S)
        \Spec(\oS^{\_phi}_{S})
            = \Spec ( \ov{\oT}^{\_phi}_S ) \bigcup (0,\dots,0)$,
        where $0$ is taken $N-|S|$ times. \\
$2)~$ Let $\_phi(x)$ be a nonnegative function.
      Then operator $\ov{\oT}^{\_phi}_{S}$ is nonnegative definite. \\
$3)~$ Let $\Gr=\Z_p$, $p$ be a prime number, and
      $S\subseteq \Gr$ be an arbitrary nonempty set.
      Suppose also that $\_phi(x)$ is a nonnegative function.
      Operator $\ov{\oT}^{\_phi}_{S}$ is positively definite iff
      there exist at least $|S|$ elements $x\in \Gr$ such that $\_phi(x) \neq 0$.
}
\label{pred:Spec_T_S}
\\
\Proof
The forth part of Proposition \ref{pred:Spec_T} implies the equality $\Spec (\oT^{\_phi}_{S}) = \Spec(\oS^{\_phi}_{S})$.
Further, the operator $\oS^{\_phi}_{S}$ has $N-|S|$ linearly  independent eigenfunctions
corresponding zero, namely $\d_{\ov{s}} (x)$, $\ov{s} \notin S$.
The last functions are
linearly independent with eigenfunctions of the restriction of the
operator  $\oS^{\_phi}_{S}$ onto $L(S)$.
Besides $\Spec ( \ov{\oS}^{\_phi}_S ) = \Spec ( \ov{\oT}^{\_phi}_S ) \subseteq \Spec (\oT^{\_phi}_{S})$
and we have proved the first part of the Proposition.

Further, let $M = \{ \_phi^{1/2} (t) e(-st) \}$, $s\in S$, $t\in \Gr$.
Then $MM^* = \ov{\oT}^{\_phi}_{S}$ and the operator $\ov{\oT}^{\_phi}_{S}$ is nonnegative definite.
Let now $\Gr=\Z_p$, $p$ be a prime number.
If there are  $x_1,\dots,x_l \in \Gr$, $l\ge |S|$ such that $\_phi(x_j) \neq 0$, $j\in [l]$,
then choose a square submatrix of $M$, say $M'$, which corresponds
to elements $x_1,\dots,x_{|S|}$.
By Theorem \ref{t:Chebotarev} $\det M' \neq 0$.
Using Binet--Cauchy's formula, we get $\det \ov{\oT}^{\_phi}_{S} \ge \prod_{j=1}^{|S|} \_phi (x_j) \cdot (\det M')^2 > 0$.
Hence all eigenvalues of the matrix $\ov{\oT}^{\_phi}_{S}$ are positive.
Remember $\ov{\oT}^{\_phi}_{S}$ is a symmetric operator, we obtain that it
is positively definite.
This completes the proof.

\Exm Let $S_1,S_2 \subseteq \Z_p$ be arbitrary sets such that $|S_1| \le |S_2|$.
Then by the last proposition, we see that the operator $\ov{\oT}^{S_2}_{S_1}$ is positively definite.
Note that  $0 \in \Spec ( \ov{\oT}^{S_1}_{S_2} )$ and $0$ has multiplicity $|S_2| - |S_1|$.
Hence the operator $\ov{\oT}^{S_1}_{S_2}$ is singular.

Formula (\ref{F:T_action}) can be rewritten for operator $\ov{\oT}^{\_phi}_{S}$  as
\begin{equation}\label{F:T_S_action}
    \langle \ov{\oT}^{\_phi}_S u, v \rangle
        =
            \sum_x (\F{\_phi^c} * u) (x) \ov{v} (x)
                =
                    \sum_x \_phi (x) \F{u} (x) \ov{\F{v} (x)} = \langle \ov{\oS}^{\_phi}_S u, v \rangle \,,
\end{equation}
where $u,v$
are
arbitrary
functions such that $\supp u, \supp v \subseteq S$.
Besides
\begin{equation}\label{F:trace}
    \tr ( \ov{\oT}^{\_phi}_{S} )
        =
            \tr (\ov{\oS}^{\_phi}_{S} )
        =
            |S| \F{\_phi} (0)
                =
                    \sum_{j=1}^{|S|} \mu_j ( \ov{\oT}^{\_phi}_{S} )
                        =
                            \sum_{j=1}^{N} \mu_j ( \oT^{\_phi}_{S} ) \,.
\end{equation}
If $\_phi$ is a real function then as was noted before
$\ov{\oT}^{\_phi}_{S}$ is a symmetric matrix.
In particular, it is a normal matrix.
Using (\ref{F:clever_convolution}),
and identities  $(S * S^c) (z) = (S * S^c) (-z)$, $\F{S^c} (\xi) = \ov{\F{S} (\xi)}$,
we get
\begin{equation}\label{F:trace_sq}
    \tr ( \ov{\oT}^{\_phi}_{S} ( \ov{\oT}^{\_phi}_{S} )^*)
        =
            \sum_z |\F{\_phi} (z)|^2 (S * S^c) (z)
                =
                    \sum_z (\_phi * \_phi^c) (z) |\F{S} (z)|^2
                =
                    \sum_{j=1}^{|S|} \mu^2_j ( \ov{\oT}^{\_phi}_{S} )
                        =
                            \sum_{j=1}^{N} \mu^2_j ( \oT^{\_phi}_{S} ) \,.
\end{equation}

The identities from the first part of Proposition \ref{pred:Spec_T} show that
there is a duality between operators $\oT^{\_phi}_{\psi}$ and $\oT^{\psi}_{\_phi}$.
Having this in mind one can suppose that
%That is why is is naturally
it should exists a dual version of operator $\ov{\oT}^{\_phi}_{S}$.
In the rest of the section we will define the dual version which we will call  $\ov{\oT}^{S}_{\_phi}$.
Also we will study eigenvalues and eigenfunctions of operators $\oT^{\_phi}_{S}$ and $\oT^{S}_{\_phi}$ in
the case $\Gr=\Z_p$, where $p$ is a prime number.

Let us make a general remark.
Let $f(x)$ be an eigenfunction of operator $\oT^{\_phi}_{\psi}$ corresponding an eigenvalue $\mu$.
In other words $f(x)$ is a non--zero function and $\oT^{\_phi}_{\psi} f = \mu f$.
Consider operator $(\oM_{\_phi} f) (x) := \_phi (x) \F{f} (x)$
and define the function $F := \oM_{\_phi} f$.
Using the third identity from (\ref{F:T_simple_op}), we get $\oT^{\psi^c}_{\_phi} F = \mu F$.
Thus, if $F(x)$ is a non--zero function then $F$ is an eigenfunction of operator $\oT^{\psi^c}_{\_phi}$
corresponding $\mu$.

\Pred
{\it
    Let $\_phi,\psi$ be arbitrary functions. \\
$1)~$     Suppose that $\oT^{\_phi}_{\psi}$ is a simple  matrix.
    Then $\oT^{\_phi^c}_{\psi^c}$ is also a simple matrix. \\
$2)~$ Let $\oT^{\_phi}_{\psi}$ be a simple matrix and $\_phi$, $\psi$ are real functions.
      Then all matrices $\oT^{\_phi^c}_{\psi}$, $\oT^{\_phi}_{\psi^c}$, $\oT^{\_phi^c}_{\psi^c}$ are simple. \\
$3)~$ Suppose that $\oT^{\_phi}_{\psi}$ is a simple matrix and $\supp \_phi = \Gr$.
    Then $\oT_{\_phi^c}^{\psi}$ is also a simple matrix.
}
\label{pred:simplicity_-1}
\\
\Proof
Let $f$ be an eigenfunction of operator $\oT^{\_phi}_{\psi}$
corresponding an eigenvalue $\mu$.
Then
\begin{equation}\label{tmp:16.09.2009_1}
    \psi(x) (\F{\_phi^c} * f) (x) = \mu f(x) \,.
\end{equation}
Apply the operator $\oC$ to the last identity and use formula (\ref{F:clever_convolution}), we get $1)$.
Let us prove the second part.
Apply the operator $\mathcal{C}$ to (\ref{tmp:16.09.2009_1})
and use formula (\ref{F:complex_op}), we obtain that $\oT^{\_phi^c}_{\psi}$ is a simple matrix.
The simplicity of another two matrices can be derived from the first part.

Further, the third identity from  (\ref{F:operators_identities}) gets $\Phi^4 = N^2 \cdot \oI$.
Hence $\Phi$ is a non--singular operator.
Let $\{ f_j \}$, $j\in [N]$ be a basis of eigenfunctions of the operator $\oT^{\_phi}_{\psi}$.
Since $\oT^{\_phi}_{\psi}$ is a simple matrix then these functions are  linearly independent.
For each $j\in [N]$ consider the function $F_j = \oM_{\_phi} f_j$.
Using the condition $\supp \_phi = \Gr$ and non--singularity of the operator $\Phi$,
we get that these functions are  linearly independent and non--zero in particular.
This completes the proof.

%Теперь разберем ситуацию группы $\Z_p$.
Now let us consider the case $\Gr = \Z_p$.

\Pred
{\it
    Let $\Gr=\Z_p$, $p$ be a prime number, and $S\subseteq \Gr$ be an arbitrary set. \\
$1)~$ Let $\psi$ be a nonnegative function.
If there are at least $|S|$ elements $x\in \Gr$ such that $\psi(x) \neq 0$ then $\oT^S_\psi$ is a simple matrix. \\
$2)~$ Let $\_phi$ be a nonnegative function.
If for all $x\in \Gr$, we have $\_phi(x) \neq 0$  then $\oT^{\_phi}_S$ is a simple matrix.
}
\label{pred:simplicity}
\\
\Proof
%Убедимся in справедливости первого пункта предложения \ref{pred:simplicity}.
Using formula (\ref{F:T}) and the identity $\F{e(ax)} = p \cdot \d_a(x)$ is is easy to see that
the functions $e(\ov{s} x)$, $\ov{s}\notin S$ are eigenfunctions of the operator $\oT^S_\psi$
corresponding to zero.
Clearly, these functions are  linearly independent.
Further, by Proposition \ref{pred:Spec_T}, we have
$\Spec (\oT^S_\psi) = \Spec (\oT^{\psi}_{S^c}) = \Spec (\oT^{\psi^c}_{S})$.
Using the third part of Proposition \ref{pred:Spec_T_S}, we see
$\oT^S_\psi$ has $|S|$ positive eigenvalues.
Hence, the correspondent eigenfunctions are linearly independent with $e(\ov{s} x)$, $\ov{s}\notin S$.
Let $f_j$, $j\in [|S|]$ be eigenfunctions of operator $\ov{\oT}^{\psi}_{S^c}$.
Clearly, $\supp f_j \subseteq S^c$, $j\in [|S|]$.
Let $F_j = M_{\psi} f_j$, $j\in [|S|]$.
We assert that these functions are linearly independent and non--zero, in particular.
Indeed, Theorem \ref{t:Tao} and condition $|\supp \psi| \ge |S^c| = |S|$
imply that the operator $M_{\psi}$ is invertible onto $L(S^c)$.
Using the argument as before, we get
$F_j$ are linearly independent functions of the operator $\oT^S_\psi$
corresponding positive eigenvalues.
Hence $\oT^S_\psi$ is a simple matrix.

Let us prove the second part of Proposition \ref{pred:simplicity}.
For each $\ov{s}\notin S$ consider the equation $M_{\_phi} f_{\ov{s}} = e(-\ov{s} x)$.
Here  $f_{\ov{s}}$ is an unknown function.
Since $\_phi(x) \neq 0$ for all $x\in \Gr$, it follows that
the equation is solvable.
Using formula (\ref{F:T}) it is easy to see that $\oT^{\_phi}_S f_{\ov{s}} = 0$ for all $\ov{s}\notin S$.
Besides, the functions $f_{\ov{s}}$, $\ov{s}\notin S$ are linearly independent.
Indeed, if $\sum_{\ov{s}\notin S} c_{\ov{s}} f_{\ov{s}} \equiv 0$
then $\sum_{\ov{s}\notin S} c_{\ov{s}} \F{f}_{\ov{s}} \equiv 0$.
Remember the definition of the functions $f_{\ov{s}}$, we have
$\sum_{\ov{s}\notin S} c_{\ov{s}} \_phi^{-1} (x) e(-\ov{s} x) \equiv 0$.
Using Theorem \ref{t:Chebotarev} for the matrix $\Phi$, we obtain that all
coefficients $c_{\ov{s}}$ equal zero.
Hence, the all the functions $f_{\ov{s}}$, $\ov{s}\notin S$ are linearly independent.
We have $|\supp \_phi| = p \ge |S|$.
By assumption $\_phi (x)$ is a non--negative function.
Using  the third part of Proposition \ref{pred:Spec_T_S}, we get all eigenvalues of the operator
$\ov{\oT}^{\_phi}_S$ are positive.
Whence, the functions $f_{\ov{s}}$, $\ov{s}\notin S$ are linearly independent with the eigenfunctions of $\ov{\oT}^{\_phi}_S$.
Thus, $\oT^{\_phi}_S$ is a simple matrix.
This completes the proof.

Finally, let us define the operator $\ov{\oT}^{S}_{\psi}$.
Let $\psi(x)$ be an arbitrary complex function.
Consider the linear space
$L^* (S) = \{ f ~:~ f = \psi (x) a(x),\, \supp \F{a} \subseteq S \}$ and, analogously,  $L^* (\ov{S})$.
Since $\oT^{S}_{\psi} f = \psi (x) (\F{ S^c \F{f^c}}) (x)$ it follows that
the space $L^* (S)$ is invariant under the action of the  operator $\oT^{S}_{\psi}$.
By definition, the restriction $\oT^{S}_{\psi}$ onto $L^* (S)$ is the operator $\ov{\oT}^{S}_{\psi}$.
Analogously, the space $L^* (S)$ is invariant under the action of the  operator $\oS^{S}_{\psi}$,
and we can define $\ov{\oS}^{S}_{\psi}$ as the restriction of the operator $\oS^{S}_{\psi}$ onto the space.
It is easy to see, that $L^* (S) = \oM_\psi L(S)$, $L^* (\ov{S}) = \oM_\psi L(\ov{S})$.
Thus, we have, in particular,
$\Spec \ov{\oT}^{S}_{\psi} = \Spec \ov{\oT}^{\psi}_{S} = \Spec \ov{\oS}^{\psi}_{S}$,
provided by the operator $M_\psi$ is invertible onto $L(S)$
(the last situation happens if $\supp \psi = \Gr$, for example, or $|\supp \psi| \ge |S|$, if $\Gr = \Z_p$).
It is not difficult to see (from the forth part of Proposition \ref{pred:Spec_T}, say)
that $\Spec \ov{\oT}^{S}_{\psi} \neq \Spec \ov{\oS}^{S}_{\psi}$ in general.

Let us define the dimension of the space $L^* (S)$.

\St
{\it
    Let $p$ be a prime number, $\Gr=\Z_p$,  $S\subseteq  \Gr$ be a nonempty set,
    and $\psi$ be an arbitrary complex function.
    Then $\dim (L^* (S)) = \min \{ |S|, |\supp \psi| \}$.
}
\\
\Proof
Since $\dim (L(S)) = |S|$ it follows that $\dim (L^* (S)) \le \min \{ |S|, |\supp \psi| \} := m$.
Let $\supp \psi = \{ a_1,\dots, a_t\}$.
Consider the case $t\le |S|$ (we can use similar arguments at the opposite situation).
Apply Corollary \ref{cor:Tao_th} to the system of linearly independent eigenfunctions of the operator $\ov{\oT}^{\psi}_S$.
To prove the inequality $\dim (L^* (S)) \ge m = t$
it is sufficiently to solve the equation
$\sum_{s\in S} c_s e(-sx) = \d_{a_i} (x)$ for each $i\in [t]$.
Here coefficients $c_s$ are unknowns.
Clearly, the last equation is solvable, because by Theorem \ref{t:Chebotarev} the matrix $\{ e(-sx) \}_{s\in S,\, x\in \supp \psi }$
has full rank.
This completes the proof.

\Exm Let $S_1,S_2 \subseteq \Z_p$ be arbitrary nonempty sets such that $|S_1| \le |S_2|$.
Then the operator $\ov{\oT}^{S_2}_{S_1}$ can be defined as the restriction of $\oT^{S_2}_{S_1}$ onto $L(S_1)$
and it can be defined as the restriction of  $\oT^{S_2}_{S_1}$ onto $L^* (S_2)$.
It is easy to see that all these definitions define the same operator.

\refstepcounter{section}
\label{Chang_gen}

\begin{center}
{\bf \arabic{section}. Some generalizations of M.--C. Chang's inequality.}
\end{center}

%Указанные простейшие
Properties of the operators $\oT^{\_phi}_\psi$, $\ov{\oT}^{\_phi}_S$,
which were considered at the previous section
allow us to prove a theorem of M.--C. Chang \cite{Ch_Fr} from the theory of large exponential sums.
Moreover, we will show that the theorem and its generalizations of different types
can be derived from Proposition \ref{st:mu_1_L_S}, which is essentially a reformulation
of Rudin's inequality  (\ref{f:Rudin}).

For a real function $\_phi(x)$ let $\_phi_0 (x) = \_phi(x)$, if $x\neq 0$
and zero otherwise.

\Pred
{\it
    Let $\Gr$ be a finite Abelian group,
    $\_phi(x)$ be a real non--zero function,
    and
    $\L \subseteq \Gr$ be a dissociated set.
    Then
    \begin{equation}\label{f:mu_1_L_S}
        |\mu_1 (\ov{\oT}^{\_phi}_{\L})| \ll \| \_phi \|_1 \l( \log ( N \| \_phi \|_{\infty} \| \_phi \|^{-1}_1 ) + 1 \r) \,.
    \end{equation}
    Let now $|\_phi (0)| = \| \_phi \|_{\infty}$.
    If $|\_phi (0)| \ge \| \_phi_0 \|_1$ then
    \begin{equation}\label{f:mu_1_L_S_second_1}
        |\mu_1 (\ov{\oT}_{\L}^{\_phi})| \ll |\_phi (0)| \log N \,.
    \end{equation}
    If $|\_phi (0)| \le \| \_phi_0 \|_1$ then
    \begin{equation}\label{f:mu_1_L_S_second_2}
        |\mu_1 (\ov{\oT}_{\L}^{\_phi})| \ll \| \_phi_0 \|_1 \l( \log ( N \| \_phi \|_{\infty} \| \_phi_0 \|^{-1}_1 ) + 1 \r) \,.
    \end{equation}
}
\label{st:mu_1_L_S}
\Proof
By assumption $\_phi(x)$ is a real function, so
$\mu_1 (\ov{\oT}_{\L}^{\_phi})$ is a real number.
Let $w$ be an arbitrary function, $\supp w \subseteq \L$, $\| w \|_2 = 1$.
Using (\ref{F:T_S_action}), we have
$$
    \sigma := \langle \ov{\oT}^{\_phi}_\L w, w \rangle = \sum_{x} |\F{w} (x)|^2 \_phi (x) \,.
$$
Let $k>0$ be an integer parameter.
By Rudin's inequality
%(see e.g. \cite{Rudin_book} or же \cite{Green_Chang2}, \cite{Sh_lacunary}, \cite{Shkr_Bourgain_group})
and H\"{o}lder's inequality, we get
\begin{equation}\label{tmp:inter_phi}
    |\sigma|^k
        \le
            \sum_x |\F{w} (x)|^{2k} \cdot \l( \sum_x |\_phi (x)|^{k/(k-1)} \r)^{k-1}
                =
                    \sigma_1 \cdot \sigma_2
                \le
\end{equation}
$$
                \le
                    C^k N \| w \|^{2k}_2 k^k \| \_phi \|_{\infty} \| \_phi \|_1^{k-1}
%                        =
                        =
                             C^k k^k \| \_phi \|_1^{k} \cdot N \| \_phi \|_{\infty} \| \_phi \|^{-1}_1 \,,
$$
where $C>0$ is the absolute constant from (\ref{f:Rudin}).
Putting $k=[\log (N \| \_phi \|_{\infty} \| \_phi \|^{-1}_1)]+1$ and using Propositions \ref{pred:Spec_T}, \ref{pred:Spec_T_S},
and Theorem \ref{t:KF_matrix}, we obtain
$$
    \mu_1 (\oT^{\L}_{\_phi}) = \mu_1 (\ov{\oT}^{\_phi}_\L) \le \sigma
        \ll
            \| \_phi \|_1 \l( \log ( N \| \_phi \|_{\infty} \| \_phi \|^{-1}_1 ) + 1 \r) \,.
$$

We need to check (\ref{f:mu_1_L_S_second_1}), (\ref{f:mu_1_L_S_second_2}).
Return to  (\ref{tmp:inter_phi}).
We have
$$
    \sigma^{1/(k-1)}_2
        =
            \sum_x |\_phi (x)|^{k/(k-1)} \le |\_phi (0)|^{k/(k-1)}
                +
                    \| \_phi_0 \|^{1/(k-1)}_\infty \| \_phi_0 \|_1 \,.
$$
By assumption $|\_phi (0)| = \| \_phi \|_{\infty}$.
If $|\_phi (0)| \ge \| \_phi_0 \|_1$ then $\sigma_2 \le 2^{k-1} |\_phi (0)|^k$.
Put $k=[\log N] + 1$ and use the arguments as before, we get (\ref{f:mu_1_L_S_second_1}).
If $|\_phi (0)| \le \| \_phi_0 \|_1$ then $\sigma_2 \le 2^{k-1} \| \_phi \|_{\infty} \| \_phi_0 \|^{k-1}_1$.
Substitute the last inequality into (\ref{tmp:inter_phi}) and put
$k = [\log ( N \| \_phi \|_{\infty} \| \_phi_0 \|^{-1}_1 )] + 1$, we obtain (\ref{f:mu_1_L_S_second_2}).
This completes the proof.

\Exm Suppose that  $\_phi(x)=S(x)$, where $S$ is an arbitrary subset of the group $\Gr$.
Then inequality (\ref{f:mu_1_L_S}) implies $|\mu_1 (\ov{\oT}^{\L}_S)| \ll |S| \log 1/\d$.

Let us obtain some applications of Rudin's inequality.
First of all derive Chang's Theorem from Proposition \ref{st:mu_1_L_S}.

\Th
{\it
    Let $\d \in (0,1]$ be a real number,
    $\L \subseteq \Gr$ be a dissociated set, and $S\subseteq \Gr$ be an arbitrary set, $|S| = \d N$.
    Then for any function $f$, $\supp f \subseteq S$, we have
    \begin{equation}\label{F:Chang_l_2_my_f}
        \sum_{x\in \L} |\F{f} (x)|^2 \ll |S| \cdot \| f \|^2_2 \log (1/\d) \,.
    \end{equation}
    In particular
    \begin{equation}\label{F:Chang_l_2_my}
        \sum_{x\in \L} |\F{S} (x)|^2 \ll |S|^2 \log (1/\d) \,.
    \end{equation}
}
\label{t:Chang_my}
%\\
\Proof
We give even two proofs.
Put $\_phi(x)$ equals $S(x)$.
Using inequality (\ref{f:mu_1_L_S}) of the last proposition, we get
$|\mu_1 (\ov{\oT}^{\L}_S)| \ll |S| \log 1/\d$.
By (\ref{F:trace}), we have $\mu_1 (\ov{\oT}^{\L}_S) > 0$.
Using Courant--Fischer's Theorem and identity (\ref{F:T_S_action}), we obtain
$$
    \sum_{x\in \L} |\F{f} (x)|^2
        =
            \langle \ov{\oT}^{\L}_S f, f \rangle
                \le
                    \mu_1 (\ov{\oT}^{\L}_S) \| f \|^2_2
        \ll
            |S| \cdot \| f \|^2_2 \log (1/\d)
$$
as required.
Let us give another proof.
Using Theorem \ref{t:KF_matrix}, formulas (\ref{F:clever_convolution}),
the second part of the Proposition \ref{pred:Spec_T} and Proposition \ref{st:mu_1_L_S}, we get
$$
    \langle \oT^{S}_{\L^c} \F{f}^c, \oT^{S}_{\L^c} \F{f}^c \rangle
        =
            \sum_{x\in \L} |(\F{S} * \F{f}) (x)|^2
                =
                    N^2 \sum_{x\in \L} |\F{f} (x)|^2
                        \le
                            \mu_1 ( \oT^{S}_{\L^c} (\oT^{S}_{\L^c})^* ) \cdot \| \F{f}^c \|^2_2
                                =
$$
$$
                                =
    N \cdot \mu_1 (\ov{\oT}^{S}_{\L^c}) \cdot N \| f \|^2_2
        \ll
            |S| \cdot \| f \|^2_2 \log (1/\d) N^2 \,.
$$
This completes the proof.

\Note We can analogously obtain a generalization of Chang's Theorem, belonging J. Bourgain \cite{Bu_new}
(see detailed  proof in \cite{Shkr_Bourgain_group}).

Our approach allows to prove new formulas.
For example, using (\ref{F_Par}),
the bound $\mu_1 (\ov{\oT}^{\L}_S) \ll |S| \log 1/\d$
and Courant--Fischer's theorem, we get
$$
    \langle (\ov{\oT}^{\L}_{S})^2 S,S \rangle
        =
            \sum_z S(z) ( \F{\L}^c * S(\F{\L}^c * S) ) (z)
                =
%                    \sum_z \ov{\F{S} (z)} \L(z) \cdot (\F{S} * (\F{S} \L)) (z)
%                        =
                            \sum_{z\in \L} \ov{\F{S} (z)} (\F{S} * (\F{S} \L)) (z)
    \le
$$
$$
    \le
        \mu^2_1 (\ov{\oT}^{\L}_{S}) \cdot \| S \|^2_2
            \ll
                |S|^3 \log^2 (1/\d) \,.
$$
Hence for any dissociated set $\L$, we have
$$
    \l| \sum_{z\in \L} \ov{\F{S} (z)} (\F{S} * \F{S} \L ) (z) \r|
            \ll
                |S|^3 \log^2 (1/\d) \,.
$$

\Note Using (\ref{F:trace}), (\ref{F:trace_sq}), we can obtain some information
about all eigenvalues $\mu_j (\ov{\oT}^{S}_{\L}) = \mu_j$, $j\in [|\L|]$,
not only $\mu_1 (\ov{\oT}^{S}_{\L})$.
First of all, by Proposition \ref{pred:Spec_T_S} all these eigenvalues
are nonnegative
(in the case $\Gr = \Z_p$, $p$ is a prime number, we have $\mu_j > 0$ for any $j$).
Secondly, from (\ref{F:trace}), we have $\sum_{j=1}^{|\L|} \mu_j = |\L| |S|$.
Thirdly, by formula (\ref{F:trace_sq}),
we get
$$
    \sum_{j=1}^{|\L|} \mu^2_j = \sum_z |\F{S} (z)|^2 (\L * \L^c) (z) = |S|^2 |\L| + \sum_{z\in \L \dotplus \L}  |\F{S} (z)|^2 \,,
$$
where $\L \dotplus \L^c = \{ \la_1 - \la_2 ~:~ \la_1,\la_2 \in \L,\, \la_1 \neq \la_2 \}$.
Finally, using an inequality of J. Bourgain (see \cite{Bu_new} or \cite{Shkr_Bourgain_group}), we obtain
$$
    \sum_{j=1}^{|\L|} (\mu_j - |S|)^2 \ll |S|^2 \log^2 (1/\d) \,.
$$
In particular, the last formula implies that
there are at most
$O(1)$ numbers $\mu_j$ such that
$\mu_j \gg |S|\log (1/\d)$.

Now let us obtain Theorems \ref{t:my_Chang_gen} and \ref{t:my_Chang_dual}.

{\bf Proof of Theorem \ref{t:my_Chang_gen}}
Let $f(x) = S(x) - \d$.
Clearly, $\F{f} (0) = 0$ and $\F{f} (x) = \F{S} (x)$, $x\neq 0$.
Define the function $F(x)$ by the formula $\F{F}(x) = |\F{S} (x)|$.
If $l=2$ then put $u(x) = (F * f) (x)$.
If $l=2k$, $k\ge 2$ then let $u(x) = (F * (f *_{k-1} f) ) * (f^c *_{k-2} f^c) (x)$.
Finally, if  $l=2k+1$, $k\ge 1$ then put $u(x) = (f*_{k} f) * (f^c *_{k-1} f^c) (x)$.
Let also $v(x) = S (x)$.
By assumption $\L$ is a dissociated set.
Hence $0\notin \L$.
Using formula (\ref{F:T_action}) and
inclusion $\supp v \subseteq S$, it is easy to see that
for any $l \ge 2$, we have
\begin{equation}\label{tmp:13.09.2009_1}
    \sigma := \langle \oT^{\L}_S u, v \rangle = \sum_{x\in \L} |\F{S} (x)|^{l+1} \,.
\end{equation}
Using the Cauchy--Schwartz inequality, the second part of Proposition \ref{pred:Spec_T}
and the bound $\mu_1 (\oT^{\L}_S) \ll |S| \log (1/\d)$, we get
$$
    \sigma^2
        \le
            \langle \oT^{\L}_S u, \oT^{\L}_S u \rangle \cdot \langle v,v \rangle
                \le
                    N \mu_1 (\oT^{\L}_S) \| u \|^2_2 \| v \|^2_2
                        \ll
                            |S|^2 \log (1/\d) \cdot \sum_{x\neq 0} |\F{S} (x)|^{2l} \,.
$$
Substitute the last identity into (\ref{tmp:13.09.2009_1}), we get
$$
    \sum_{x \in \L} |\F{S} (x)|^{l+1}
        \ll
            |S| \log^{1/2} (1/\d) \cdot \l( \sum_{x\neq 0} |\F{S} (x)|^{2l} \r)^{1/2}
$$
as required.
This completes the proof.

\Note
Let $\d,\a \in (0,1]$ be real parameters, $\a \le \d$, and $\L \subseteq \Z_N$ be an arbitrary dissociated set,
such that $|\L| \gg (\d/\a)^2 \log (1/\d)$.
Developing the approach from papers \cite{Green_Chang_exact,Green_Chang1}, in article \cite{Sh_exp2} (see Theorem 2.8)
the following set $S$ was constructed.
Suppose that there are some restrictions onto parameters $\d,\a$,
and $|\L| \ll (\d/\a)^2 \log (1/\d)$.
Then exists a set $S$,  $S\subseteq \Z_N$, $|S| = \d N$  such that  \\
$1)~$ For all $x\neq 0$, we have $|\F{S}(x)| \ll \a N$. \\
$2)~$ For each $x\in \L \bigsqcup (-\L)$, the following holds $|\F{S} (x)| \gg \a N$. \\
$3)~$ For any $x\neq 0$, $x\in \L \bigsqcup (-\L)$, we have $|\F{S} (x)| \ll \eps N$, where $\eps = \a^2 \d^{-1}$. \\
The set $S$ gives us an example showing that inequality (\ref{f:my_Chang_gen}) of Theorem \ref{t:my_Chang_gen} is sharp.
Indeed,
by $2)$ and the inequality $|\L| \gg (\d/\a)^2 \log (1/\d)$, we obtain
\begin{equation}\label{tmp:13.09.2009_I}
    \sum_{x \in \L} |\F{S} (x)|^{l+1}
        \gg
            |\L| (\a N)^{l+1}
                \gg
                    (\a N)^l \d^2 \a^{-1} \log (1/\d) \,.
\end{equation}
From the other hand, using Parseval's identity and $1)$, $3)$, we get
$$
    \sum_{x\neq 0} |\F{S} (x)|^{2l}
        =
            \sum_{x\in \L \bigsqcup (-\L)} |\F{S} (x)|^{2l}
                +
            \sum_{x\notin \L \bigsqcup (-\L),\, x\neq 0} |\F{S} (x)|^{2l}
                \ll
$$
$$
                \ll
                    |\L| (\a N)^{2l} + (\eps N)^{2l-2} \d N^2
                        \ll
                            |\L| (\a N)^{2l} \,.
$$
To prove the last inequality we have supposed that $\a \ll \d \cdot \d^{1/(2l-2)}$.
Thus
\begin{equation}\label{tmp:13.09.2009_II}
    |S| \log^{1/2} (1/\d) \cdot \l( \sum_{x\neq 0} |\F{S} (x)|^{2l} \r)^{1/2}
        \ll
            \d N \log^{1/2} (1/\d) |\L|^{1/2} (\a N)^l
                \ll
                    (\a N)^l \d^2 \a^{-1} \log (1/\d)
\end{equation}
and we see that assuming the condition  $\a \ll \d \cdot \d^{1/(2l-2)}$, we have
lower bound (\ref{tmp:13.09.2009_I}) is coincident (up to constants) to upper bound (\ref{tmp:13.09.2009_II}).

Now let us prove Theorem \ref{t:my_Chang_dual}.
We obtain even a tiny stronger result.

\Th
{\it
    Let $\Gr$ be a finite Abelian group, and $l$ be a positive integer, $l\ge 2$.
    Let also $\L \subseteq \Gr$ be dissociated set,
    and $S_1,\dots, S_l \subseteq \Gr$ be arbitrary sets.
    Then
\begin{equation}\label{f:my_Chang_dual'}
    \sum_{x\in \L} (S_1 * S_2 * \dots * S_l)^2 (x)
        \ll
            \frac{|S_l|}{N} \l( \sum_{x} \prod_{j=1}^{l-1} | \F{S}_j (x) |^2 \r) \cdot \log N \,.
\end{equation}
    If
    $\prod_{j=1}^{l-1} |S_j|^2 \ge \sum_{x\neq 0} \prod_{j=1}^{l-1} |\F{S}_j (x)|^2$
 then
\begin{equation}\label{f:my_Chang_dual'_1}
    \sum_{x\in \L} (S_1 * S_2 * \dots * S_l)^2 (x)
        \ll
            \frac{|S_l|}{N} \l( \prod_{j=1}^{l-1} |S_j|^2 \r) \cdot \log N \,.
\end{equation}
    If
    $\prod_{j=1}^{l-1} |S_j|^2 \le \sum_{x\neq 0} \prod_{j=1}^{l-1} |\F{S}_j (x)|^2$
    then
\begin{equation}\label{f:my_Chang_dual'_2}
    \sum_{x\in \L} (S_1 * S_2 * \dots * S_l)^2 (x)
        \ll
            \frac{|S_l|}{N} \l( \sum_{x\neq 0} \prod_{j=1}^{l-1} |\F{S}_j (x)|^2 \r) \cdot \log N \,.
\end{equation}
Suppose that $S_1=S_2=\dots=S_{l-1}=S$, $|S| \le N/2$ and
$|S|^{2l-2} \le \sum_{x\neq 0} |\F{S} (x)|^{2l-2}$.
Then
\begin{equation}\label{f:my_Chang_dual'_4_1}
    \sum_{x\in \L} ( (S *_{l-2} S) * S_l)^2 (x)
        \ll
            \frac{l |S_l|}{N} \l( \sum_{x\neq 0} |\F{S} (x)|^{2l-2} \r) \cdot \log |S| \,.
\end{equation}
    Finally, always
\begin{equation}\label{f:my_Chang_dual'_simple}
    \sum_{x\in \L} (S_1 * S_2)^2 (x)
        \ll
            |S_1| |S_2| \cdot \log ( \min\{ |S_1|, |S_2| \} ) \,.
\end{equation}
}
\label{t:my_Chang_dual'}
{\bf Proof of Theorem \ref{t:my_Chang_dual'}}
Let $f(x) = S_l (x)$, $\_phi(x) = \prod_{j=1}^{l-1} \F{S}^c_j (x)$.
Using the arguments similar to the second variant of the proof of Theorem \ref{t:Chang_my}, we get
$$
    \langle \oT^{\_phi}_{\L^c} f^c, \oT^{\_phi}_{\L^c} f^c \rangle
        =
            \sum_{x\in \L} |(\F{\_phi} * f) (x) |^2
                =
                    N^2 \sum_{x\in \L} (S_1 * S_2 * \dots * S_l)^2 (x)
                        \le
                            \mu_1 ( \oT^{\_phi}_{\L^c} (\oT^{\_phi}_{\L^c})^* ) \cdot \| f^c \|^2_2
                                =
$$
$$
                                =
    N \cdot \mu_1 (\ov{\oT}^{|\_phi|^2}_{\L^c}) \cdot \| f \|^2_2
        \ll
            N |S_l| \mu_1 (\ov{\oT}^{\_phi^*}_{\L^c}) \,,
$$
where $\_phi^* (x) = |\_phi (x)|^2 = \prod_{j=1}^{l-1} |\F{S}^c_j (x)|^2$.
We have $\| \_phi^* \|_{\infty} = \_phi^* (0) = \prod_{j=1}^{l-1} |S_j|^2$
and $\| \_phi^* \|_1 \ge \prod_{j=1}^{l-1} |S_j|^2$.
Proposition \ref{st:mu_1_L_S} gives us some inequalities for the quantity $\mu_1 (\ov{\oT}^{\_phi^*}_{\L^c})$.
Applying (\ref{f:mu_1_L_S}), (\ref{f:mu_1_L_S_second_1}),
we get (\ref{f:my_Chang_dual'}) and (\ref{f:my_Chang_dual'_1}), correspondingly.
If $\prod_{j=1}^{l-1} |S_j|^2 \le \sum_{x\neq 0} \prod_{j=1}^{l-1} |\F{S}_j (x)|^2$ then
$\| \_phi^*_0 \|_1 \ge |\_phi^* (0)| = \| \_phi^* \|_{\infty}$ and
inequality  (\ref{f:mu_1_L_S_second_2}) implies (\ref{f:my_Chang_dual'_2}).
In the case $l=2$ the quantity $\| \_phi^* \|_1$ can be computed.
Indeed, by Parseval's identity, we obtain $\| \_phi^* \|_1 = |S_1| N$.
Whence,  inequality (\ref{f:my_Chang_dual'_simple}) holds.
Finally, we need to check (\ref{f:my_Chang_dual'_4_1}).
By assumption $|S|^{2l-2} \le \sum_{x\neq 0} |\F{S} (x)|^{2l-2}$.
Thus we can use inequality (\ref{f:mu_1_L_S_second_2}) of Proposition \ref{st:mu_1_L_S}
to estimate $\mu_1 (\ov{\oT}^{\_phi^*}_{\L^c})$.
By assumption $|S| \le N/2$.
The last bound and Parseval's identity imply $\sum_{x\neq 0} |\F{S} (x)|^2 = N|S| - |S|^2 \ge 2^{-1} N|S|$.
Using H\"{o}lder's inequality, we obtain $\| \_phi^*_0 \|_1 \ge 2^{-(l-1)} N |S|^{l-1}$
and inequality
%This implies
(\ref{f:my_Chang_dual'_4_1}).
This concludes  the proof.

\Cor
{\it
    Let $\Gr$ be a finite Abelian  group, $r$ be a positive integer.
    Let also $\L \subseteq \Gr$ be a dissociated set,
    and $S_1,S_2 \subseteq \Gr$ be arbitrary sets.
    Suppose that for any $x\in \L$ the following holds $(S_1 * S_2) (x) \ge r$.
    Then
\begin{equation}\label{f:large_convolution}
    |\L| \ll r^{-2} |S_1| |S_2| \cdot \log ( \min\{ |S_1|, |S_2| \} )  \,.
\end{equation}
}
\label{cor:large_convolution}

\Note Let $\Gr = (\Z / p\Z)^n$, $p$ be a prime number
(usefulness of considering such groups was discussed in survey \cite{Green_finite_fields}),
and $S_1=S_2 = P$, where $P$ be a linear subspace of $\Gr$.
Clearly, any dissociated set $\L \subseteq P$ has the cardinality at most $\log |P|$
and there are dissociated sets $\L \subseteq P$ such that $|\L| \gg \log |P|$.
From the other hand from (\ref{f:my_Chang_dual'_simple}), we have $r=|P|$ and
$\sum_{x\in \L} (S_1 * S_2)^2 (x) = |\L| |P|^2 \ll |P|^2 \log |P|$.
Hence $|\L| \ll \log |P|$.
Thus, at least in the situation when the parameter $r$ is large
inequality (\ref{f:my_Chang_dual'_simple})
and consequently Corollary \ref{cor:large_convolution}
are sufficiently sharp.

If $r$ is a small number then there is a different bound for the cardinality of $\L$.
We thank to S. Yekhanin and K. Talwar for pointed us the fact.
For simplicity, let $\Gr = (\Z / 2\Z)^n$, and $S_1=S_2=S$.
Take $p=c r^{-1/2} \log^{1/2} |S|$, where $c>0$ is an absolute constant, and
suppose that $r\gg \log |S|$.
Let us choose a random subset $S' \subseteq S$ such that
any element $x$  belongs to $S'$ with probability $p$.
Clearly, the expectation of the cardinality of the set  $S'$ equals $p|S|$.
By assumption  for any $x\in \L$, we have $(S*S) (x) \ge r$.
Hence choosing the constant  $c$,  we get $\L \subseteq S' + S'$ with positive probability
(probability of the event that $\L$ not in $S' + S'$ does not exceed $|\L| (1-p^2)^r \le |S|^2 (1-p^2)^r$).
It is easy to see (the proof was also suggested to us by S. Yekhanin and can be found e.g. in \cite{Konyagin_Sh})
that $|\L| \ll |S'|$.
Hence we have with positive probability that $|\L| \ll |S'| \ll p |S| \ll |S| r^{-1/2} \log^{1/2} |S|$.
If $r$ is small then the last bound is better than (\ref{f:large_convolution}).

\end{document}